\input amstex
\input amsppt.sty
\magnification=\magstep1
\advance\vsize-1cm\voffset=-0.5cm\advance\hsize1cm\hoffset0cm
\NoBlackBoxes
\def\R{\Bbb R}   

\topmatter
\title A characterization of submanifolds by a homogeneity condition \endtitle
\author Arkadiy Skopenkov \endauthor
\address
Department of Differential Geometry, Faculty of Mechanics and
Mathematics, Moscow State University, 119992, Moscow, Russia,
and Independent University of Moscow, B. Vlasyevskiy, 11, 119002,
Moscow, Russia.
e-mail: skopenko\@mccme.ru \endaddress
\abstract 
A very short proof of the following smooth homogeneity theorem of D. Repovs, 
E. V. Scepin and the author is presented. 
Let $N$ be a locally compact subset of a smooth manifold $M$. 
Assume that for each two points $x,y\in N$ there exist their neighborhoods $Ux$ 
and $Uy$ in $M$ and a diffeomorphism $h:Ux\to Uy$ such that $h(x)=y$ and 
$h(Ux\cap N)=Uy\cap N$. Then $N$ is a smooth submanifold of $M$. 
\endabstract
\endtopmatter

\document
Which shape could have a sheath so that it would be possible to draw a sabre 
out of it? 
Mathematical formulation of this question leads to the following notion. 
A subset $N$ of the 3-dimensional (or $m$-dimensional Euclidean) space is 
called  {\it Riemannian ambient homogeneous} if for each two points $x,y\in N$ 
there exists an isometry $h:\R^3\to\R^3$ that maps $x$ to $y$ and $N$ to $N$. 
It is well known that {\it each Riemannian ambient homogeneous curve in the
3-dimensional space is either a straight line or a circle or a spiral line.} 

Which shape could have a cable so that it would be possible to draw a wire out 
of it (it is allowed to flex a wire but not to break it)?
Mathematical formulation of this question leads to the following notion. 
A subset $N$ of a smooth manifold $M$ is called {\it smoothly ambient 
homogeneous} if {\it for each two points $x,y\in N$ there exists a 
diffeomorphism $h:M\to M$ that maps $x$ to $y$ and $N$ to $N$.}
(The Theorem and Applications below is interesting and non-trivial even for the 
case when $M=\R^m$ or even when $M$ is the plane.) 
In this paper 'smooth' means 'differentiable'; see though remarks at the end. 

Any smooth submanifold of a smooth manifold (in particular, a graph of a 
smooth function $\R\to\R$) is smoothly ambient homogeneous.
In this note the converse is proved. 
The proof is simpler than in [RSS93, RSS96, RS00] (although it uses the same 
ideas).

\smallskip
{\bf Theorem.} 
{\it Let $N$ be a locally compact subset of a smooth manifold $M$. 
If $N$ is smoothly ambient homogeneous then $N$ is a smooth submanifold of 
$M$.}

\smallskip
{\bf Applications.} 
(1) If a graph of a continuous function $\R^{m-1}\to\R$ is smoothly ambient 
homogeneous, then the function is smooth. 
(A function having infinite derivative at some point is differentiable at 
this point.)

(2) The Cantor set cannot be smoothly ambient homogeneously embedded into 
$\R^m$. 


(3) It is known that manifolds are homogeneous and that converse is false (the 
Cantor set is a counterexample). 
The Theorem shows that the property of being a {\it submanifold} is
equivalent to the {\it ambient homogeneity} property. Cf. [Gl68].

(4) Using the Theorem it is convenient to prove that some groups are Lie 
groups. 
E.g. it implies the Cartan theorem stating that {\it any closed subgroup of a 
Lie group is a Lie group.}

(5) For applications to the Hilbert-Smith conjecture see [RSS96, RSS97]. 


(6) The Theorem allows to reduces of the following result [MZ55, Theorem 3 on 
p. 208-209] to its simpler case $m=1$ (i.e. to the Cauchy equation 
$h(s+t)=h(s)+h(t)$): {\it if a one-parameter group $\{h^t\}_{t\in \R}$ of
diffeomorphisms of an $m$-dimensional manifold depends continuously
on a parameter $t$, then it depends smoothly on the parameter.} (V.
I. Arnold suggested the result as a problem in 1980's.)

\smallskip
{\it Proof of the Theorem.}
The property of being a smooth submanifold in $M$ is a local one. 
Hence we may assume that $M=\R^m$. 

Let $B^{m+1}_l$ be the interior of an $m$-ball of radius $1/l$ without the 
center. 
Let $B^0_l:=\emptyset$. 
For $1\le k\le m$ let $B^k_l$ the open cone over the $(1/l^2)$-neighborhood of 
$(k-1)$-hemisphere in the $(m-1)$-sphere of radius $1/l$: 
$$B^k_l\quad:=\quad\{(x_1,\dots,x_m)\in\R^m\ |\ -l^2x_k<|x|<1/l\quad\text{and}
\quad l^2|x_i|<|x|\text{ for }k<i\le m\},$$
where $|x|=\sqrt{x_1^2+\dots+x_m^2}$. 
Denote by $O_m$ the group of orthogonal transformations of $\R^m$. 
Take the greatest $k\ge0$ such that 
$$\text{(*)\quad for each $x\in N$ there exist $l$ and $A\in O_m$ such that 
$(x+AB^k_l)\cap N=\emptyset$}.$$ 
(Informally this means that $N$ is '$(m-k)$-dimensionally Lipschitz'.)
Such $k$ exists because (*) is true for $k=0$. 
If $k=m+1$, then $N$ consists of isolated points and the Theorem is proved. 
Hence we may assume that $k\le m$. 

Take a sequence $\{A_l\}$ everywhere dense in $O_m$.
Denote
$$B_l:=A_lB^k_l\quad\text{and}\quad N_l:
=\{x\in N\ |\ (x+B_l)\cap N=\emptyset\}.$$
By (*) we have $N=\cup_{l=1}^\infty N_l$.  
It is easy to check that {\it $N_l$ is closed in $N$} (see the details in 
[RSS96, Lemma 3.1]).
Therefore by the Baire Category Theorem some $N_l$ contains a non-empty 
open in $N$ set.

So there exists a point $x\in N$ and an $m$-dimensional cube $I^m$ of diameter 
less than $1/l$ and with the centre $x$, for which $N':=N\cap I^m\subset N_l$.
Then
$$(**)\quad
[(y+B_l)\cup(y-B_l)]\cap N'=\emptyset\quad\text{for each}\quad y\in N'.$$
Indeed, if $z\in(y-B_l)\cap N'$, then $y\in(z+B_l)\cap N'$, which is impossible.

Since $N$ is locally compact, we may assume that $N'$ is compact. 
Also we may assume that $I^m=I^{m-k}\times I^k$ so that $I^k$ is parallel 
to the $k$-dimensional plane $A_l(\R^k\times\vec0)$. 
Let $p:I^m\to I^{m-k}$ be the projection. 

\smallskip
{\it First case: $p(N')$ contains open in $I^{m-k}$ set.} 
(This is automatically true for $k=m$, when everything is already clear by 
(**).)   
We may assume that this set is $I^{m-k}$ itself (by changing $I^m$ to the 
product of a part of this set with $I^k$). 
From (**) it follows that $N'$ is a graph of a Lipschitz map $q:I^{m-k}\to I^k$.
Therefore $q$ has a differentiability point [Fe69, Theorem 3.1.6].
Now the smooth ambient homogeneity implies that $q$ is differentiable. 
Hence by smooth ambient homogeneity $N$ is smooth submanifold.

\smallskip
{\it Second case: $p(N')$ does not contain any open in $I^{m-k}$ set.} 
(Then $k<m$.)   
It follows that there exists a point $a\in I^{m-k}-p(N')$ close enough to the 
centre of $I^{m-k}$.
Since $p(N')$ is compact, the distance from $a$ to $p(N')$ is non-zero and 
there is a point $z\in N'$ such that $|a-p(z)|$ equals to this distance.  
Then the open ball $D\subset I^{m-k}$ with the centre $a$ and radius $|a-p(z)|$ 
does not intersect $p(N')$.
So $p^{-1}(D)\cap N'=\emptyset$.
Clearly, 
$$(z+B_l)\cup(z-B_l)\cup p^{-1}(D)\supset z+A_lB^{k+1}_s\quad\text{for some }s.
$$ 
This and (**) imply that $(z+A_lB^{k+1}_s)\cap N=\emptyset$.
Since $N$ is smoothly ambient homogeneous, it follows that for each $x\in N$ 
there exists 
a diffeomorphism $h:\R^m\to\R^m$ mapping $z$ to $x$ and $N$ to $N$. 
Then 
$$h(N\cap(z+A_lB^{k+1}_s))\supset x+AB^{k+1}_u\quad\text{for some}
\quad A\in O_m.$$ 
Hence (*) remains true if we replace $k$ by $k+1$. 
This contradicts to the maximality of $k$. \qed 

\smallskip
Note that the Cantor set {\it can} be {\it continuously} or {\it Lipschitz} 
ambiently homogeneously embedded into the plane and $\R^m$ [MR99], see also 
[DRS89, DR95].
Hence the analogues of the Theorem for continuous or Lipschitz categories are 
false. 

The analogue of the Theorem (and of Application (1)) is true for 
$C^1$-category. 
The proof is analogous. 
(At the end of the first case it should be additionally noted that the 
derivative of a differentiable mapping has a continuity point. 
Then ambient $C^1$-homogeneity implies that $N$ is a $C^1$-submanifold.) 
I conjecture that {\it the analogue of the Theorem holds for $C^r$-category 
when $r\ge2$ and for analytic category.}
In spite of [RSS96, RSS97], I do not have a proof of this conjecture for 
$C^r$-category when $r\ge2$. 

The Theorem (with an analogous proof) remains valid if by an ambient 
differentiable homogeneity we understand the following propety: 
{\it for each two points $x,y\in N$ there exist their 
neighborhoods $Ux$ and $Uy$ in $M$ and a diffeomorphism $h:Ux\to Uy$ that maps 
$x$ to $y$ and $Ux\cap N$ to $Uy\cap N$.}

I am grateful to A. Efimov for useful discussions. 

\bigskip
\centerline{\bf References.} 

\smallskip
\eightpoint
[DRS89] D. Dimovski, D. Repov\v{s} and E. V.\v S\v cepin, 
$C^\infty$-homogeneous closed curves on orientable closed surfaces, 
Geometry and Topology, ed. G. M. Rassles and G. M. Stratopoulos, 
1989 World Scientific Publ. Co Singapore, pp. 100-104.

[DR95] D. Dimovski and D. Repov\v{s}, On homogeneity of compacta in manifolds,  
Atti. Sem. Mat. Fis. Univ. Modena, XLIII (1995) 25--31. 

[Fe69] H. Federer, Geometric Measure Theory, Springer, Berlin, 1969. 

[Gl68] H.~Gluck, Geometric characterisation of differentiable manifolds in 
Euclidean space, II, Michigan Math.~J. {\bf 15:1} (1968), 33--50. 

[MR99] J.Male\v{s}i\v{c} and D. Repov\v{s},
On characterization of Lipschitz manifolds,
New Developments In Differential Geometry, J. Szenthe, Ed., Kluwer, 
Dordrecht 1999, pp. 265-277.

[MZ55] D. Montgomery and L. Zippin, Topological Transformation Groups, 
Princeton, Princeton Univ. Press, 1955. 

[RS97] D. Repov\v{s} and E. V. Shchepin, A proof of the Hilbert-Smith conjecture 
for actions by Lipschitz maps, Math. Ann. {\bf 308} 1997, 361--364.

[RS00] E. V. Shchepin and D. Repov\v{s}, On smoothness of compacta.
Jour. of Math. Sci., {\bf 100(6)} 2000, 2716--2726.

[RSS93] D.~Repov\v{s}, A.~B.~Skopenkov and E.~V.~\v{S}\v{c}epin, 
A characterization of $C^1$-homogeneous subsets of the plane, Boll. Unione 
Mat. Ital., {\bf 7-A} (1993), 437--444. 

[RSS96] D. Repov\v{s}, A. B. Skopenkov and E. V.\v S\v cepin,
$C^1$-homogeneous compacta in $\R^n$ are $C^1$-submanifolds of $\R^n$,
Proc. Amer. Math. Soc. {\bf 124:4} (1996), 1219--1226.

[RSS97] D.~Repov\v s, A.~B.~Skopenkov and E.~V.~\v S\v cepin, Group actions on 
manifolds and smooth ambient homogeneity, Jour. of Math. Sci. (New York), 
{\bf 83:4} (1997), 546--549. 

\enddocument 
\end